\documentclass[11pt,twoside]{amsart}
\usepackage{amssymb}

\newtheorem{thm}{Theorem}[section]

\newtheorem{cor}[thm]{Corollary}

\newtheorem{defn}[thm]{Definition}

\newtheorem{example}[thm]{Example}
\newtheorem{rem}[thm]{Remark}
\newtheorem{conj}[thm]{Conjecture}

\newcommand{\skipit}[1]{{}}
\newcommand{\prfend}{\hbox to7pt{\hfil}
\par\vskip-\baselineskip\hbox to\hsize
{\hfil\vbox {\hrule width6pt height6pt}}\vskip\baselineskip}

\newcommand{\myarrow}[2]{\hbox to #1pt{\hfil$\to$\hfil}{\hskip-#1pt{\raise
10pt\hbox to#1pt{\hfil$\scriptscriptstyle #2$\hfil}}}}

\begin{document}

\title{An improved Multiplicity Conjecture for codimension three Gorenstein algebras}

\author{Juan C. Migliore
\and
Uwe Nagel
\and
Fabrizio Zanello}

\address{
Department of Mathematics,
University of Notre Dame,
Notre Dame, IN 46556}
\address{
Department of Mathematics, University of Kentucky, 715 Patterson Office Tower, Lexington, KY 40506-0027}
\address{
Department of Mathematics, Royal Institute of Technology, S-100 44 Stockholm, Sweden}

\maketitle

\markboth{J.\ Migliore, U.\ Nagel, F.\ Zanello}{An improved multiplicity conjecture}

\begin{abstract}
The Multiplicity  Conjecture is a deep problem relating the multiplicity (or degree) of a Cohen-Macaulay standard graded algebra with certain extremal graded Betti numbers in its minimal free resolution.  In the case of level algebras of codimension three, Zanello has proposed a stronger conjecture.  We prove this conjecture in the Gorenstein case.
\end{abstract}

\section{Introduction}\label{intro}
The third author has recently proposed an improvement of  the so-called Multiplicity Conjecture, for the case of codimension three level algebras.  The purpose of this note is to prove the first case of this improved conjecture, namely that of codimension three Gorenstein  algebras.  Furthermore, we characterize the cases when the improved conjecture for Gorenstein algebras is sharp.

Let $R$ be a polynomial ring over a field.  Consider a Cohen-Macaulay standard graded algebra, $R/I$, of codimension $r$ with minimal free resolution
\[
0 \rightarrow \bigoplus_{j \geq 0} R(-j)^{\beta_{r,j}} \rightarrow \bigoplus_{j \geq 0} R(-j)^{\beta_{r-1,j}} \rightarrow \dots \rightarrow \bigoplus_{j \geq 0} R(-j)^{\beta_{1,j}} \rightarrow R \rightarrow R/I \rightarrow 0.
\]
Let $M_i = \max \{ j \ | \ \beta_{i,j} \neq 0 \}$ and let $m_i = \min \{ j \ | \ \beta_{i,j} \neq 0 \}$.  The Multiplicity Conjecture of Herzog, Huneke and Srinivasan says that
\[
\frac{1}{r!} \prod_{i=1}^r m_i \leq e(R/I) \leq \frac{1}{r!} \prod_{i=1}^r M_i,
\]
where $e(R/I)$ is the multiplicity of $R/I$.
It was shown by Huneke and Miller \cite{HM} that if the minimal free resolution is {\em pure} (i.e.\ for all $i$, $\beta_{i,j} \cdot \beta_{i,k} = 0$ whenever $j \neq k$) then both bounds above are true (and thus are sharp).  This conjecture has been shown to hold in a number of important cases, and it has been extended beyond the Cohen-Macaulay situation.  We refer to the expository paper \cite{FS} for a more detailed (but at this point no longer up-to-date) history of the progress on the problem.  However, of greatest interest for us here is the following result of the first and second author together with T.\ R\"omer:

\begin{thm}[\cite{MNR1}, Theorem 1.4] \label{ht 3 gor bd}
Let $R/I$ be a standard graded Gorenstein algebra of codimension three.
Then the following lower and upper bounds hold:
\begin{enumerate}\label{eq-gor}
\item
$e(R/I)\geq \frac{1}{6} m_1 m_2 m_3 + \frac{1}{6}(M_3-M_2)^2
  (M_2-m_2+M_1-m_1)$;
\item
$e(R/I) \leq \frac{1}{6} M_1 M_2 M_3 - \frac{1}{12} M_3(M_2-m_2+M_1-m_1)$.
\end{enumerate}
\end{thm}

As an immediate consequence we get not only the fact the Multiplicity Conjecture holds for codimension three Gorenstein ideals, but in fact
that either bound of the Multiplicity Conjecture  is sharp
for codimension three Gorenstein ideals if and only if the resolution is pure:

\begin{cor}[\cite{MNR1}, Corollary 1.5] \label{cor-sharp-gor}
Let $R/I$ be a graded Gorenstein algebra of codimension three.
Then the following
conditions are equivalent:
\begin{enumerate}
\item $e(R/I) = \frac{1}{6} m_1 m_2 m_3$;
\item $e(R/I) = \frac{1}{6} M_1 M_2 M_3 $;
\item $R/I$ has a pure minimal free resolution.
\end{enumerate}
\end{cor}

Because of this result, it was conjectured in  \cite{MNR2} and in
\cite{HZ} that in all cases, the bounds of the Multiplicity
Conjecture are sharp if and only if the resolution is pure.  Special
cases were proven in \cite{MNR2}, \cite{HZ} and also in \cite{MR}.

Corollary \ref{cor-sharp-gor} of course gives a strengthening of the original conjecture.  Other variants have been suggested as well.  One such idea was carried out by Francisco \cite{francisco} and  then from a different point of view by the third author \cite{zanello}.    We give a brief outline of this approach.  Suppose that $R/I$ is a Cohen-Macaulay standard graded algebra, and consider a minimal free resolution of $R/I$ as above.  If there is a redundant term in this minimal free resolution (i.e.\ a copy of $R(-j)$ occurring in consecutive free modules), then formally removing this term from both modules gives a new ``Betti sequence'' that may or may not actually occur as a resolution for some Cohen-Macaulay algebra $R/I'$.  Nevertheless, one can formally compute the multiplicity as if this sequence were exact, and it will equal the multiplicity of the original algebra $R/I$.  However, the $m_i$ potentially increase, and the $M_i$ potentially decrease.  If the multiplicity formally computed with the new ``Betti sequence" satisfies the bounds with the new Betti numbers, then $e(R/I)$ also satisfies the bounds of the Multiplicity Conjecture.  Thus the strengthened conjecture is that obtained by replacing the $m_i$ and the $M_i$ by these new Betti numbers.

One drawback to this approach, in general,  is that there may well
be choices in the cancellation of terms.  For instance, there may be
three consecutive copies of the same $R(-j)$, so choosing one pair
gives a different result than choosing the other pair.  Francisco's
approach is to make the convention that the rightmost cancellation
takes priority.

Zanello's approach differs in two ways.  First,  he restricts to the
case of {\em level} algebras of codimension 3.  Recall that a level
algebra of codimension 3 is a Cohen-Macaulay ring for which
$\beta_{3,j} = 0$ for all but one value of $j$.  With this
restriction, there is no choice in the cancellation.  Second, he
notes that in this case the strengthened conjecture can be phrased
entirely in terms of the Hilbert function of $R/I$.  We give the
details, notation and precise statement of the strengthened
conjecture in the next section. However, we would like to point out
that in computing the ``Betti numbers" for the improved bound, we
may well be looking at a ``Betti diagram" that does not exist.  We
will see that this can happen even in the Gorenstein case.  It was
shown by Francisco that in his approach, the reduction to diagrams
that do not actually exist can produce {\em numerical}
counterexamples to the Multiplicity Conjecture.  That this does not
happen in the Gorenstein case (as we show) or more generally in the
level case (conjecturally) is remarkable.

The main result of this paper (Theorem \ref{main theorem})  is that
if $R/I$ is a codimension 3 Gorenstein graded algebra then it
satisfies the strengthened conjecture.  We use results of Diesel as
well as the fact that both bounds of the original Multiplicity
Conjecture are known to hold (\cite{HS} Theorem 2.4,  \cite{MNR1}
Theorem 1.4).  Furthermore, we characterize those Gorenstein graded
algebras for which the strengthened conjecture is sharp.


\section{Hilbert functions and cancelling redundant terms}

In the recent paper \cite{zanello}, the third author  sought to find
new conjectural bounds  which are stronger than those of the
Multiplicity Conjecture (in general), and which could be expressed
purely in terms of the Hilbert function of $R/I$.  Because of the
nature of this approach, the conjecture is stated only for
codimension three level algebras.  However, as we will see, a
similar  idea was proposed by Francisco also in higher codimension.

Let $A = R/I$ be a codimension three Artinian algebra with Hilbert function
\[
\underline{h} =  (1,3,h_2,\dots,h_{c-1}, h_c).
\]
  It is a standard fact (see \cite{FL}) that the third difference of $\underline{h}$ has the following connection with the graded Betti numbers of $A$:
\[
\Delta^3 \underline{h}(t) :=  h_t - 3 h_{t-1} + 3 h_{t-2} - h_{t-3} = \beta_{2,t} - \beta_{1,t} - \beta_{3,t}
\]
for all $t \geq 1$.  Notice that under most  circumstances, this
formula does not guarantee that we know the values of any of the
Betti numbers.  This formula does not pick up the existence of
``ghost" or redundant terms in the minimal free resolution (i.e.\
terms $R(-t)$ that are repeated in the second and either the first
or the third free module (or both) in the resolution).

However, we can draw the following conclusions:

\begin{enumerate}
\item if $\Delta^3 \underline{h}(t) > 0$  then $\beta_{2,t} > 0$.

\item The initial degree, $m_1$, of $I$ can be read from the Hilbert function -- it is
\[
m_1 = \min \left  \{ i \ | \ h_i < \binom{i+2}{2} \right \}.
\]
In this degree of course $\beta_{2,t} = \beta_{3,t} =  0$, so
$\Delta^3 \underline{h} < 0$.  In fact $m_1 = \min \{ t \ | \
\Delta^3 \underline{h}(t) < 0 \}$ and $\beta_{1,m_1} = - \Delta^3
\underline{h} (m_1)$.

\item Since $A$ is Artinian, we have
\[
-\Delta^3 \underline{h} (c+3) = \beta_{3,c+3} = h_c \ \ \hbox{ and } \ \ \Delta^3 \underline{h}(t) = 0 \ \hbox{ for all } \ t > c+3.
\]
  If we further assume that $A$ is {\em level}, then $\beta_{3,t} = 0$ for all $t$ except $t = c+3$. Hence, in this case, for all $t < c+3$ we have
\[
\Delta^3 \underline{h}(t) <0 \ \ \Rightarrow \ \ \beta_{1,t} > 0.
\]
\end{enumerate}

\begin{defn}  Assume that $A$ is level.
We set

\begin{itemize}
\item $n_1 := \min  \{ t \ | \ \Delta^3 \underline{h} (t) < 0  \}$. \\

\item $n_2 := \min \{ t >0 \ | \ \Delta^3 \underline{h} (t) > 0 \}$. \\

\item $N_1 := \max \{ t \leq c+1 \ | \ \Delta^3 \underline{h}(t) < 0 \}$. \\

\item $N_2  :=  \max \{ t \ | \ \Delta^3 \underline{h}(t) >0 \}$.

\end{itemize}
\end{defn}

\begin{rem} \label{note1} {\rm
\begin{enumerate}
\item Whether or not $A$ is level, we have $n_1 = m_1$.

\item Since $A$ is level, we have $N_2 = M_2$.

\item We also have $n_2 \geq m_2$ and $N_1 \leq M_1$.  What might prevent equality is the existence of redundant terms.

\item If $A$ is level then, in the Multiplicity Conjecture, we have $m_3 = M_3 = c+3$, which can be read uniquely from the Hilbert function.
\end{enumerate}}
\end{rem}

Because of Remark \ref{note1}, the third author made the following conjecture:

\begin{conj}[\cite{zanello}, Conjecture 3.1]  \label{zanello conjecture}
If $R/I$ is level of codimension three then
\[
\frac{1}{3!} n_1 n_2 (c+3) \leq e(R/I) \leq \frac{1}{3!} N_1 N_2 (c+3).
\]
\end{conj}
Note that this is a strengthening of the Multiplicity Conjecture.
As pointed out in the introduction, the idea is similar to one
described by C.\ Francisco in \cite{francisco}.


\section{A proof of the conjecture for codimension three Gorenstein algebras}

In this section we prove the main result of this paper:

\begin{thm} \label{main theorem}
Conjecture \ref{zanello conjecture} holds for codimension three Gorenstein algebras.
\end{thm}

\begin{proof}  We first need some notation and background results.

Using the structure theorem of Buchsbaum and Eisenbud \cite{BE}, Diesel \cite{diesel} gave a complete description of the degree matrices that can occur in the middle of the minimal free resolution of a codimension three Gorenstein algebra.  (She gave a number of consequences of this result, which do not concern us here.)  We set our notation, mostly following \cite{diesel}.

Let $A = R/I$ be a graded, codimension 3 Artinian Gorenstein algebra over $R := k[x,y,z]$, where $k$ is a field.  From \cite{BE} we know that a minimal generating set for $I$ has an odd number of elements.  So suppose that $R/I$ has a minimal free resolution
\[
0 \rightarrow R(-s) \rightarrow \bigoplus_{i=1}^n R(-p_i) \stackrel{B}{\longrightarrow} \bigoplus_{i=1}^n R(-q_i) \rightarrow R \rightarrow R/I \rightarrow 0
\]
where $n$ is odd.  Buchsbaum and Eisenbud showed that $B$ can be chosen to be skew-symmetric, and the minimal generators are the pfaffians of $B$.  The socle of $R/I$ occurs in degree $s-3$, the $q_i$ are the degrees of the minimal generators of $I$ and the $p_i$ are the degrees of the minimal relations among the generators.  We make the convention that
\[
\begin{array}{l}
q_1 \leq q_2 \leq \dots \leq q_n \\ \\
p_1 \geq p_2 \geq \dots \geq p_n.
\end{array}
\]
Notice that $p_i + q_i = s$ for all $i = 1,\dots,n$.  We consider the degree matrix
\[
\delta B =
\left [
\begin{array}{cccccccccccccccccc}
p_1 - q_1 & p_2 - q_1 & \dots & p_n - q_1 \\
p_1 - q_2 & p_2 - q_2 & \dots & p_n - p_2 \\
p_1 - q_3 & p_2 - q_3 & \dots & p_n - q_3 \\
\vdots & \vdots & & \vdots \\
p_1 - q_n & p_2 - q_n & \dots & p_n - q_n
\end{array}
\right ]
\]
The entries of $\delta B$ are increasing as one moves up and to the left.

The key fact for us is summarized by Diesel in \cite{diesel}, first paragraph of Section 2.3.  We paraphrase her observations as follows:

\begin{quotation}
{\em
Under our hypotheses and notation, if $\delta B$ contains a {\em non-diagonal}  entry that is zero then there is a second such entry.  In this case we can allow these entries to be non-zero constants $c$ and $-c$.  The resulting ideal of pfaffians will still be Gorenstein, of height 3, and minimally generated by $n-2$ elements.  }
\end{quotation}

This says that whenever $p_i = q_j$ for $i \neq j$ then the corresponding (numerically redundant) terms in the resolution have associated to them a second ``dual'' pair of numerically redundant terms, and these four summands can be removed, with the strong conclusion that the result is still the minimal free resolution of a Gorenstein algebra of codimension three.

We will consider the lower bound of Conjecture \ref{zanello conjecture}.  The upper bound is proved in a completely analogous fashion.

Given the graded Betti numbers of $R/I$, the computation of $n_2$ is performed as follows.  Start with $p_n$ (the smallest syzygy degree).  If there is any $q_i$ equal to $p_n$, with $i \neq n$, then we will also have $p_i = q_n$.  Thus the summands in the minimal free resolution corresponding to $q_i, p_n, p_i, q_n$ can be removed, and what remains is still the set of graded Betti numbers of a codimension 3 Gorenstein algebra $A' = R/I'$.  If no such $q_i$ exists then $n_2 = p_n$.  Continuing in this manner, we can end in one of two ways:

\begin{enumerate}

\item at some point we have computed $n_2$ by the nonexistence of $q_i$, or

\item we come to a point where the only possible cancellation comes from $q_n$ (suitably reindexed).
\end{enumerate}

In the first case, note that $n_2 = m_2 (A')$ for this new Gorenstein algebra $A'$ obtained by formal cancellation.  But we know that the Multiplicity Conjecture holds for $A'$, by Theorem  \ref{ht 3 gor bd}.  So
\[
\begin{array}{rcl}
\frac{1}{3!} n_1  n_2 (c+3) & = & \frac{1}{3!} m_1 m_2(A') m_3 \\
& \leq & e(A') \\
& = & e(A).
\end{array}
\]
Hence Conjecture \ref{zanello conjecture} holds for $A = R/I$.

Now suppose that (2) holds above.  So the graded Betti numbers allow some numerical cancellation that does not occur for any Gorenstein algebra.  However, diagonal entries in $\delta B$ correspond to
pairs $p_i, q_i$, while under our assumption we are always working with the syzygy of least degree.  Hence we have for the last Gorenstein algebra constructed, that $p_n = q_n$; that is, the syzygy of least degree coincides with the generator of largest degree.  Furthermore, clearly we have reduced the problem  to a situation where $\beta_{2,p_n} = \beta_{1,q_n} = 1$ in order for no non-diagonal cancellations to be possible in these degrees.
We conclude that the Betti diagram  must have the form

\begin{center}
\begin{tabular}{r|cccccccccccccccccc}
 & 1 & $n$ & $n$ & 1 \\ \hline
 & $1$ & - & - & - \\
&& \vdots & \vdots \\
& - & $a_1$ & - & - \\
&& \dots & \dots & \\
& - & $a_{r-1}$ & - & - \\
& - & $a_r$ & 1 & - \\
& - & 1 & $a_r$ & - \\
& - & - & $a_{r-1}$ & - \\
& & \vdots & \vdots \\
& - & - & $a_1$ & - \\
& - & \vdots & \vdots \\
 & - & - & - & 1
\end{tabular}

\end{center}

This Betti diagram exists (thanks to Diesel), but the computation of
$n_2$ requires one more cancellation, even though such an algebra
does not exist (since it has to be Gorenstein, but has an even
number of minimal generators).  However, note that the Betti diagram
above is   quasi-pure, and the diagram obtained by removing the 1's
is even pure.  It was shown by Herzog and Srinivasan for
Cohen-Macaulay algebras with quasi-pure resolution (\cite{HS},
Theorem 1.2), and by the first two authors and R\"omer  for
Cohen-Macaulay modules of rank 0 with quasi-pure resolution
(\cite{MNR2}, Theorem 4.2), that the Multiplicity Conjecture holds
in those cases.  Furthermore, as noted by Francisco
(\cite{francisco}, immediately after Theorem 4.7), the proof of this
result ``is numerical: Any potential quasi-pure resolution below
that of a lexicographic ideal satisfies the bounds; there is no need
for there to exist a module with that resolution."

We conclude that the lower bound of Conjecture \ref{zanello
conjecture} holds.  As remarked earlier, the upper bound is proved
similarly, and in fact the steps are dual to those given for the
lower bound (we start with the generator of largest degree).
\end{proof}

As one might expect, Conjecture \ref{zanello conjecture} can be sharp even for Gorenstein algebras that do not have pure resolution (in contrast to the situation of Corollary \ref{cor-sharp-gor}).  It depends only on the Hilbert function, and on whether the diagram obtained after performing all possible numerical cancellation is pure.  We omit the proof, which is immediate from the proof of Theorem \ref{main theorem}; the only additional ingredient is Corollary \ref{cor-sharp-gor}.

\begin{cor} \label{sharp}
For a Gorenstein algebra of codimension three, the following are equivalent:

\begin{enumerate}
\item $e(R/I) = \frac{1}{3!} n_1 n_2(c+3)$.

\item $e(R/I) = \frac{1}{3!}N_1 N_2(c+3)$.

\item $n_1 = N_1$ and $n_2 = N_2$.  (That is, the Betti diagram reduces to a pure one via formal cancellation.)

\item For all $1 \leq t \leq c+2$, there is exactly one $t_1$ and one $t_2$ for which $\Delta^3 \underline{h}(t_1) < 0$ and $\Delta^3 \underline{h}(t_2) > 0$ (all other $\Delta \underline{h}(t) = 0$); furthermore, $-\Delta^3 \underline{h}(t_1) = \Delta^3 \underline{h}(t_2)$ and $t_1 + t_2 = c+3$.

\end{enumerate}
\end{cor}

\begin{example} {\rm
As remarked, the condition in Corollary \ref{sharp} does not quite say that the minimal free resolution of the Gorenstein algebra has to be pure.  For instance, the Hilbert function
\[
1, \ \ 3, \ \ 6, \ \ 6, \ \ 3, \ \ 1
\]
leads to a third difference
\[
1, \ \ 0, \ \ 0, \ \ -4, \ \ 0, \ \ 4, \ \ 0, \ \ 0, \ \ -1,
\]
from which we immediately see that $n_1 = N_1 = 3$ and $n_2 = N_2 = 5$, and Conjecture \ref{zanello conjecture} is sharp.  However, the general minimal free resolution for such a Hilbert function is
\[
0 \rightarrow R(-8) \rightarrow
\begin{array}{c}
R(-5)^4 \\
\oplus \\
R(-4)
\end{array}
\rightarrow
\begin{array}{c}
R(-4) \\
\oplus \\
R(-3)^4
\end{array}
\rightarrow R \rightarrow R/I \rightarrow 0
\]
which is clearly not pure.}
\end{example}

\begin{rem}{\rm
One can check that the example above illustrates that there are situations when both bounds of Conjecture \ref{zanello conjecture} are sharp while neither bound of Theorem \ref{ht 3 gor bd} is sharp.  On the other hand, any Gorenstein graded algebra for which the minimal free resolution is not pure and does not have redundant terms clearly gives an instance where Theorem \ref{ht 3 gor bd} is sharper than Conjecture~\ref{zanello conjecture}.  }

\end{rem}

{\ }\\
\\
\\
{\bf Acknowledgements.} The research contained in this paper was performed during the third author's visit to the first author at the University of Notre Dame, and that visit was supported by a grant of the Vetensk\aa psradet (Swedish Research Council) and a grant of the Department of Mathematics of the University of Notre Dame. Moreover, the third author was funded by the G\"oran Gustafsson Foundation.\\
\\


\begin{thebibliography}{ll}


\bibitem{BE} D.\ Buchsbaum and D.\ Eisenbud, {\em Algebra Structures for Finite
Free Resolutions, and some Structure Theorems for Ideals of Codimension 3},
Amer.\ J.\ of Math.\ {\bf 99} (1977), 447--485.

\bibitem{diesel} S.~Diesel, {\em Irreducibility and Dimension Theorems for
Families of Height 3 Gorenstein Algebras}, Pacific J.\ Math. {\bf 172}
(1996), no. 2, 365--397.

\bibitem{francisco} C.\ Francisco,  {\em New approaches to bounding the multiplicity of an ideal},
J.\ Algebra  {\bf 299} (2006),  no.\ 1, 309--328.

\bibitem{FS} C.\ Francisco and H.\ Srinivasan, {\em Multiplicity Conjectures}.  Available for download at \newline  {\tt www.math.missouri.edu/~chrisf/multiplicity-final.pdf}.

\bibitem{FL} R. Fr\"oberg and D.\ Laksov, {\em Compressed Algebras}, Conference on Complete Intersections in Acireale, Lecture Notes in Mathematics, No.\ 1092 (1984), 121--151, Springer-Verlag.

\bibitem{HS} J.\ Herzog and H.\ Srinivasan, {\em Bounds for multiplicities},
Trans.\ Amer.\ Math.\ Soc. {\bf 350} (1998), no.\ 7, 2879--2902.

\bibitem{HZ}
J.\ Herzog  and Zheng, {\em Notes on the multiplicity conjecture},
Collectanea Math.\ {\bf 57}, 2 (2006), 211--226.

\bibitem{HM}
C.\ Huneke and M.\ Miller,
{\em A note on the multiplicity of Cohen-Macaulay algebras with pure
  resolutions}.
Can.\ J.\ Math.\ {\bf 37} (1985), 1149--1162.

\bibitem{MNR1} J.\ Migliore, U.\ Nagel and T.\ R\"omer, {\em The multiplicity
conjecture in low codimensions},  Math.\ Res.\ Lett.\ {\bf 12} (2005), 731-748.

\bibitem{MNR2} J.\ Migliore, U.\ Nagel and T.\ R\"omer, {\em Extensions of the Multiplicity Conjecture}, to appear in Trans.\ Amer.\ Math.\ Soc.  Available for download at {\tt http://front.math.ucdavis.edu/math.AC/0505229.}

\bibitem{MR} R.\ Mir\'o-Roig, {\em A note on the multiplicity of determinantal ideals},
to appear in J.\ Algebra.  Available for download at {\tt http://front.math.ucdavis.edu/math.AC/0504077}.

\bibitem{zanello} F.\ Zanello, {\em Improving the bounds of the Multiplicity Conjecture:
the codimension 3 level case}, J.\ Pure Appl.\ Algebra {\bf 209}
(2007), no.\ 1, 79--89.


\end{thebibliography}
\end{document}